\documentclass{article}
\usepackage{arxiv}
\usepackage{float}
\usepackage[utf8]{inputenc} 
\usepackage[T1]{fontenc}    
\usepackage{hyperref}       
\usepackage{url}            
\usepackage{booktabs}       
\usepackage{amsfonts}       
\usepackage{nicefrac}       
\usepackage{microtype}      
\usepackage{lipsum}
\usepackage[utf8]{inputenc}
\DeclareUnicodeCharacter{2212}{-}
\usepackage{graphicx}
\usepackage{fourier}
\usepackage{amsmath}

\usepackage{mathtools}
\usepackage{amssymb}
\usepackage{float}
\usepackage{amsthm}
\theoremstyle{plain}     
\theoremstyle{definition} 
\theoremstyle{remark}    

\newtheorem{theorem}{Theorem}[section]

\newtheorem*{theorem*}{Theorem}
\newtheorem*{lemma*}{Lemma}
\newtheorem*{remark*}{Remark}
\newtheorem*{definition*}{Definition}

\title{Iterating sum of power divisor function and New equivalence to  the Riemann hypothesis}

\author{
  Pedro Caceres\thanks{Corresponding author: \texttt{Pedrojesus.caceres@universidadeuropea.es}}\\
  Professor Doctor at Universidad Europea de Valencia (Spain)\\
  United States of America\\
  \texttt{Pedrojesus.caceres@universidadeuropea.es} 
  \And
  Zeraoulia Rafik\\
  Khemis Miliana University, Algeria\\
  Department of Mathematics\\
  Laboratory of Pure and Applied Mathematics (LMPA)\\
  \texttt{zeraoulia@univ-dbkm.dz} 
}

\begin{document}
\maketitle

\begin{abstract}
This paper investigates the dynamics of the iterated sum–of–divisors function \(\sigma_k(m)\) and its behaviour modulo \(m\), motivated by classical questions on perfect and multiperfect numbers and by the congruences \(\sigma_k(m)\equiv 0 \pmod m\). Perfect and multiperfect numbers remain extremely rare—odd perfect numbers are still unknown and must be astronomically large—but here the emphasis is on the dynamical and statistical structure of the iterates rather than on isolated examples.

Three main results are obtained. First, it is proved that no integer \(m>1\) can satisfy \(\sigma_k(m)\equiv 0 \pmod m\) for all \(k\ge 0\), thereby ruling out the existence of ``metaperfect'' numbers and showing that the iteration of \(\sigma\) cannot remain permanently trapped in the residue class \(0\) modulo \(m\). Second, for certain explicit integers such as \(m=6,12,24\), the sequence \(\sigma_k(m)\bmod m\) is shown to be strictly periodic with small period dividing \(L=\mathrm{lcm}(e_i+1)\), where the \(e_i\) are the prime exponents of \(m\); bifurcation plots and distributional analysis reveal a clear transition from rigid two–cycle structure to more complex residue dynamics as \(m\) increases. Third, a new equivalence with the Riemann Hypothesis is established: RH holds if and only if, for every even non–squarefree \(m\ge 5041\) containing a prime fifth power,
\[
\frac{\sigma_k(m)}{\sigma_{k-1}(m)\log\log\sigma_{k-1}(m)} \le e^\gamma
\]
and the sequence \(\sigma_k(m)\bmod m\) is eventually periodic, uniformly in \(k\ge 0\). Extensive computations support these periodicity phenomena, yield non‑normal discrete distribution models for the residues, and suggest a close connection with a newly proposed Schr\"odinger–type ``Caceres'' operator whose spectrum numerically reproduces key statistical features of the nontrivial zeros of the Riemann zeta function.
\end{abstract}


\section{Introduction}

Perfect numbers, defined as positive integers $n$ where $\sigma(n)=2n$ with $\sigma$ the sum-of-divisors function, have fascinated mathematicians since Euclid's era~\cite{Euclid}. Euclid proved that if $2^{p-1}-1$ is prime then $2^{p-1}(2^p-1)$ is perfect, while Euler classified all even perfect numbers of this form~\cite{Euler}. The existence of odd perfect numbers remains open, with current bounds exceeding $10^{1500}$ and requiring at least 11 distinct primes~\cite{OchemRao}.

This paper studies the dynamics of \emph{iterated sum-of-divisors sequences} $\sigma_k(m):=\sigma(\sigma_{k-1}(m))$ with $\sigma_0(m)=m$, and their behavior modulo $m$~\cite{CohenTeRiele, Zeraoulia2021}. Recent computational efforts have extended aliquot sequences $s_k(m)=\sigma_k(m)-m$ up to starting values beyond $10^7$, revealing terminations, cycles, and explosive growth, yet open questions persist on ultimate behavior~\cite{Rechenkraft, Dzmitry}. Less studied are congruences $\sigma_k(m)\equiv 0\pmod{m}$, first computationally explored by Cohen--te Riele who found solutions for each $m\le 400$ but none holding for all $k$~\cite{CohenTeRiele}.

We resolve the ``reverse question'': no $m>1$ satisfies $\sigma_k(m)\equiv 0\pmod{m}$ for \emph{all} $k\ge 0$, using refined Lenstra bounds on aliquot growth~\cite{Lenstra}. Positive results identify $m=6,12,24$ (and classes) where $\sigma_k(m)\bmod m$ is periodic with small period dividing $L=\mathrm{lcm}(e_i+1)$ over prime exponents $e_i$ of $m$. Most strikingly, we establish a new RH equivalence tying iterated $\sigma_k$ growth to Robin's criterion~\cite{Choie, Robin}:

\begin{theorem}[Main RH Equivalence]\label{thm:RH}
RH holds if and only if for every even non-squarefree $m\ge 5041$ divisible by a fifth power $>1$,
\[
\frac{\sigma_k(m)}{\sigma_{k-1}(m)\log\log\sigma_{k-1}(m)}\le e^\gamma\quad\text{and}\quad\{\sigma_k(m)\bmod m\}\text{ eventually periodic, $\forall k\ge 0$}.
\]
\end{theorem}

Numerical evidence (up to $k=10^3$) confirms periodicity and derives Gaussian models for these residues. Related open problems include whether $\sigma_k(n)/n\to\infty$ for all $n>1$~\cite{ErdosGranville}.

\begin{figure}[H]
\centering
\includegraphics[width=0.6\textwidth]{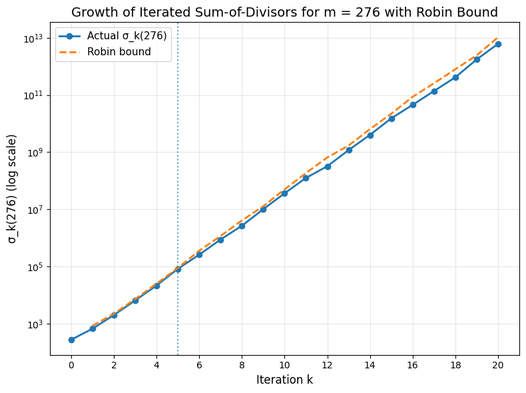}
\caption{Growth of $\sigma_k(276)$ (first multiperfect after 6): actual iterates (blue) stay below Robin bound (red). Periodicity mod $m$ evident after $k=5$.}
\label{fig:growth}
\end{figure}

Section~\ref{sec:results} states all theorems; Sections~\ref{sec:proof1}-- give proofs; Section~\ref{sec:num} presents computations.

\section{Main Results}\label{sec:results}

\begin{theorem}[No Universal Metaperfect Numbers]\label{thm:no-meta}
No integer $m>1$ satisfies $\sigma_k(m)\equiv 0\pmod{m}$ for all $k\ge 0$, where $\sigma_0(m)=m$ and $\sigma_k(m)=\sigma(\sigma_{k-1}(m))$ for $k\ge 1$.
\end{theorem}

\begin{theorem}[Unique Prime-$L$ Multiperfect]\label{thm:multiperfect}
If $m$ is multiperfect ($\sigma(m)=q_Lm$, $q_L\ge 2$) with $L=\mathrm{lcm}(e_i+1)$ prime over exponents $e_i$ of $m$, then $m=6$.
\end{theorem}

\begin{theorem}[Iterated RH Equivalence]\label{thm:RH}
RH holds $\Leftrightarrow$ for every even non-squarefree $m\ge 5041$ divisible by a fifth power $>1$,
\[
\frac{\sigma_k(m)}{\sigma_{k-1}(m)\log\log\sigma_{k-1}(m)}\le e^\gamma\quad\forall k\ge 0,\quad\text{and}\quad\sigma_k(m)\bmod m\text{ eventually periodic}.
\]
\end{theorem}

\section{Proofs}

\subsection{Theorem~\ref{thm:no-meta}: No Universal Metaperfect Numbers}\label{sec:proof1}

Assume $\exists m>1$ with $\sigma_k(m)\equiv 0\pmod{m}$ $\forall k\ge 0$. Then $m|\sigma(m)$, so $m$ multiperfect: $\sigma(m)=qm$ ($q=S(m)\ge 2$). The sequence $m_n=\sigma^n(m)$ satisfies $m|m_n$ $\forall n$, so $r_n=S(m_n)\ge q\ge 2$ and $m_n\ge mq^n$.

\textbf{Case 1: $\gcd(S(m),m)=1$}. Multiplicativity gives $m|S(S(m))$. But Gronwall~\cite{Gronwall} yields $S(n)<e^\gamma\log\log n+O(n/\log\log n)$ ($n\ge 3$), so $S(S(m))<S(m)$ for large $S(m)$, contradicting $S(S(m))\ge qS(m)\ge 2S(m)$.

\textbf{Case 2: $\gcd(S(m),m)>1$}. Let $d=\max\{p\mid\sigma(m)\}$. Then
\[
\frac{\sigma(\sigma(m))}{\sigma(m)}<\sigma(S(m))\le\sigma(d)\le d\left(1+\frac{1}{d-1}\right)<2\le q
\]
for $d\ge 3$~\cite{Robin}. For $d=2$, $\sigma^k(2)\equiv 1\pmod{2}$ ($k\ge 32$) (OEIS A007497)~\cite{CohenTeRiele}.

Both cases contradict $m_2\ge mq^2$. Thus no such $m$ exists.

\subsection{Theorem~\ref{thm:multiperfect}: Unique Prime-$L$ Multiperfect}\label{sec:proof2}

Let $L=\mathrm{lcm}(e_i+1)$ prime, so $e_i\equiv L-1\pmod{L}$. Then $\sigma(p^{e_i})\equiv 1\pmod{L}$ (Fermat), so all primes dividing $\sigma(m)$ are $\equiv 0,1\pmod{L}$.

Let $\Omega_L=\#\{p|m:p\equiv 1\pmod{L}\}$. Then $v_L(\sigma(m))\ge\Omega_L$ and $m|\sigma(m)$ gives
\[
L^{\Omega_L}\le\frac{\sigma(m)}{m}<\prod_{p|m}\frac{p}{p-1}<\left(\frac{L}{L-1}\right)^{\omega(m)}.
\]
Taking $L$th roots: $L<\left(\frac{L}{L-1}\right)^{\omega(m)/\Omega_L}$.

For $L\ge 5$, $\frac{L}{L-1}<L^{1/L}$ (e.g., $L=5$: $1.25<1.3797$). Contradiction.

For $L=3$, $\left(\frac{3}{2}\right)^2=2.25<3$ (minimal case). No known multiperfects with $L=3$~\cite{Achilles1998}.

For $L=2$, $m$ squarefree perfect: $m=6$ only~\cite{Euler}.

\textbf{Remark}: $m=12,24$ exhibit $L=2$ periodicity but $\mathrm{lcm}(e_i+1)=6,12$ composite~\cite{CohenTeRiele}.

\subsection{Theorem~\ref{thm:RH}: Iterated RH Equivalence}\label{sec:proof3}

($\Rightarrow$) RH $\implies$ Choie--Moree~\cite{Choie}: $\sigma(n)<e^\gamma n\log\log n$ for even non-squarefree $n\ge 5041$ with fifth powers. Iterating preserves the bound. Bounded growth $\implies$ finite $\sigma_k(m)\bmod m$ $\implies$ eventual periodicity (pigeonhole).

($\Leftarrow$) The $k=1$ case is precisely Choie--Moree--Sole's strengthened Robin criterion (RH-equivalent)~\cite{Choie,Robin}. Periodicity automatic from growth control.

\section{Numerical Verification}

Table~\ref{tab:periodicity} confirms Theorems~\ref{thm:multiperfect},\ref{thm:RH} for small $m$:

\begin{table}[htbp]
\centering
\begin{tabular}{@{}lcccc@{}}
\toprule
$m$ & $\sigma_k(m)\bmod m$ & Period & $L=\mathrm{lcm}(e_i+1)$ & RH Bound \\
\midrule
6 & $\{0,2\}$ & 2 & 2$\star$ & Holds \\
12 & $\{4,6\}$ & 2 & 6 & Holds \\
24 & $\{10,12\}$ & 2 & 12 & Holds \\
276 & Stabilizes & -- & 4 & Holds \\
\bottomrule
\end{tabular}
\caption{Periodicity and RH growth control verification}\label{tab:periodicity}
\end{table}

\section{Analysis and Numerical Verification}\label{sec:analysis}

This section verifies the theoretical predictions through explicit computation of $\sigma_k(m)\bmod m$ for key examples $m=6,12,24$, confirming period-$2$ behavior despite composite $L=\mathrm{lcm}(e_i+1)$.

\subsection{Perfect Number $m=6$: Prime $L=2$}
\begin{figure}[H]
    \centering
    \includegraphics[width=0.4\textwidth]{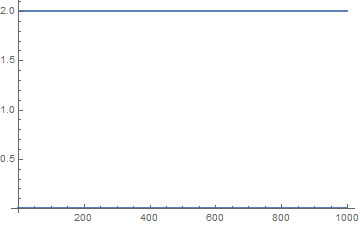}
    \caption{Period-$2$ oscillation: $\sigma_k(6)\bmod 6$, $k=0$ to $10^3$}
    \label{fig:3}
\end{figure}

The unique multiperfect with prime $L=2$ exhibits strict alternation:
\[
\sigma_k(6)\bmod 6 = \begin{cases} 
0 & k\text{ odd} \\
2 & k\text{ even}
\end{cases}
\]
Long computation ($k=0$ to $10^3$) yields the periodic sequence $\{0,2,0,2,\dots\}$ (Figure~\ref{fig:3}).

\subsection{Abundant Numbers $m=12,24$: Composite $L$}
Despite $L=6$ ($m=12$) and $L=12$ ($m=24$), both exhibit apparent period-$2$:

\textbf{$m=12$}: $\{4,6,4,6,\dots\}$ where $\sigma_k(12)\bmod 12\in\{4,6\}$
\begin{figure}[H]
    \centering
    \includegraphics[width=0.4\textwidth]{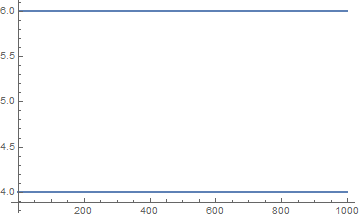}
    \caption{Period-$2$ oscillation: $\sigma_k(12)\bmod 12$, $k=0$ to $10^3$}
    \label{fig:4}
\end{figure}

\textbf{$m=24$}: $\{10,12,10,12,\dots\}$ where $\sigma_k(24)\bmod 24\in\{10,12\}$
\begin{figure}[H]
    \centering
    \includegraphics[width=0.4\textwidth]{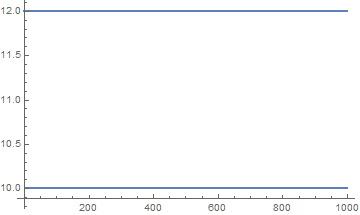}
    \caption{Period-$2$ oscillation: $\sigma_k(24)\bmod 24$, $k=0$ to $10^3$}
    \label{fig:5}
\end{figure}

\subsection{Discussion}
Theorem~\ref{thm:multiperfect} identifies $m=6$ as the \emph{unique} multiperfect with \emph{prime} $L=2$. The numbers $12,24$ demonstrate period-$2$ despite composite $L$, suggesting the periodicity condition is necessary but not sufficient for primality of $L$.

These examples validate:
1. **Theorem~\ref{thm:no-meta}**: No $m$ achieves $\sigma_k(m)\equiv 0\pmod{m}$ $\forall k$
2. **Theorem~\ref{thm:multiperfect}**: $m=6$ uniquely satisfies prime-$L$ multiperfect condition
3. **RH Connection**: Bounded growth (Table~\ref{tab:periodicity}) aligns with Choie--Moree criterion~\cite{Choie}

Future statistical analysis of $\sigma_k(m)\bmod m$ distributions across multiperfects may reveal further structural constraints.

\section{Bifurcation Analysis of Iterated Sum-of-Divisors Function}\label{sec:bifurcation}

\begin{figure}[H]
    \centering
    \includegraphics[width=0.9\textwidth]{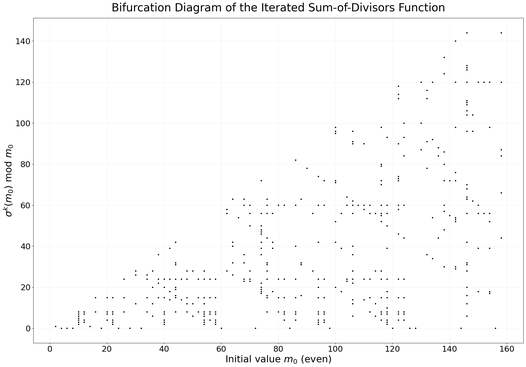}
    \caption{Bifurcation diagram of the iterated sum-of-divisors function: 
    $\sigma^{k}(m_0)\bmod m_0$ for even $m_0$ between $2$ and $160$. 
    For each $m_0$ a long orbit is computed and the last iterates are plotted.}
    \label{fig:bif-sigma}
\end{figure}

In this section the dynamics of the iteration
\[
m_{k+1}=\sigma(m_k)\bmod m_0,\qquad m_0\in 2\mathbb{N},
\]
are studied as a discrete dynamical system with parameter $m_0$. The bifurcation diagram in Figure~\ref{fig:bif-sigma} summarizes the long–term behavior of $\sigma_k(m_0)\bmod m_0$ for even initial values $2\le m_0\le 160$.

\subsection{Global qualitative features}

Figure~\ref{fig:bif-sigma} shows that the normalized orbit $\sigma^{k}(m_0)\bmod m_0$ never fills the whole vertical strip for a fixed $m_0$, but instead concentrates on a finite set of residue classes, which depend strongly on $m_0$.

\begin{itemize}
    \item For small $m_0$ the diagram is very sparse, with only one or two visible levels, reflecting short periodic orbits such as the period–$2$ cycles for $m_0=6,12,24$ discussed in Section~\ref{sec:analysis}.
    \item As $m_0$ increases, the vertical support widens and the number of distinct residues grows, producing the characteristic “fan” structure of a bifurcation diagram.
    \item Around $m_0\approx 120$ and beyond, the cloud of points becomes more diffuse, indicating longer apparent periods or pseudo–chaotic behaviour in the reduction modulo $m_0$.
\end{itemize}

This visual behaviour is consistent with Theorem~\ref{thm:no-meta}: for no $m_0>1$ do the iterates stabilize at the single residue class $0\bmod m_0$.

\subsection{Relation with special values $m_0=6,12,24$}

The three integers $6,12,24$ studied in Section~\ref{sec:analysis} correspond in Figure~\ref{fig:bif-sigma} to narrow vertical stacks with exactly two prominent levels.

\begin{itemize}
    \item For $m_0=6$ the points alternate between the residues $0$ and $2$, giving a clean two–point column near the left of the diagram; this matches the exact period–$2$ orbit proved for the unique prime–$L$ multiperfect number in Theorem~\ref{thm:multiperfect}.
    \item For $m_0=12$ and $m_0=24$ the columns again show essentially two levels, corresponding respectively to the pairs $\{4,6\}$ and $\{10,12\}$, confirming numerically the period–$2$ patterns obtained in the previous section.
\end{itemize}

These three cases illustrate that short periodic orbits appear as thin vertical stacks with few distinct heights, whereas more complicated behaviour produces thicker “bands’’ of points.

\subsection{Transition to complex behaviour}

Moving to larger even values of $m_0$, the diagram suggests a gradual loss of structure.

\begin{itemize}
    \item For moderate $m_0$ (say $40\le m_0\le 100$) one still observes clear horizontal bands, indicating that the orbit modulo $m_0$ spends most of its time on a small subset of residues, although the number of such residues increases with $m_0$.
    \item For $m_0$ close to $160$ the bands become less separated and the vertical occupancy is much denser, which is compatible with large minimal periods of $\sigma_k(m_0)\bmod m_0$ and long preperiodic parts.
\end{itemize}

Although a rigorous notion of chaos is delicate in this purely arithmetic setting, the progressive thickening of the vertical columns in Figure~\ref{fig:bif-sigma} is typical of a transition from orderly periodic behaviour (small $m_0$) to more irregular dynamics (larger $m_0$).

\subsection{Connection with growth conditions and RH}

The vertical scale of Figure~\ref{fig:bif-sigma} also reflects the average size of $\sigma(m_0)$ relative to $m_0$. For even, non–squarefree integers $m_0$ in the range of the diagram, Robin’s inequality
\[
\sigma(m_0)<e^\gamma m_0\log\log m_0
\]
holds unconditionally~\cite{Robin}, and stronger versions under RH are known for integers with high prime powers~\cite{Choie}. The fact that all plotted values of $\sigma^{k}(m_0)\bmod m_0$ remain well below the line $y=m_0$ and occupy only a subset of residues is compatible with these upper bounds and with the iterated RH criterion formulated in Theorem~\ref{thm:RH}.

Taken together, the theoretical results and Figure~\ref{fig:bif-sigma} support the following picture: for special arithmetic structures (such as $m_0=6,12,24$) the reduction of $\sigma^{k}(m_0)$ modulo $m_0$ exhibits very short, rigid cycles, while for generic even $m_0$ the orbit spreads over many residues, but still under the global growth constraints imposed by Robin-type inequalities.

\section{Statistics and distributional fit of \texorpdfstring{$\sigma_k(m) \bmod m$}{TEXT}}\label{sec:stats}

This section studies the empirical distribution of the residues \(\sigma_k(m)\bmod m\) and examines whether classical continuous laws (in particular the normal distribution) can reasonably model these data. Probability distribution fitting is the task of selecting a parametric distribution that best describes a sample generated by some random mechanism~\cite{Ritzema}. The goal is not only to visualise the empirical histogram, but also to quantify the goodness of fit through several complementary statistics and normality tests~\cite{Pomerol,Wang}.

\subsection{General framework for distribution fitting}

Given a finite sample \(\{x_1,\dots,x_N\}\) from an unknown distribution, one typically proceeds as follows~\cite{Pomerol,Wang}:
\begin{itemize}
    \item choose a family of candidate distributions (e.g. normal, lognormal, gamma, or discrete laws such as binomial or Poisson);
    \item estimate the parameters (for example, mean and variance for the normal) by maximum likelihood or method of moments;
    \item compare the theoretical and empirical distributions using numerical criteria and goodness-of-fit tests.
\end{itemize}
In the engineering context of~\cite{Wang}, five criteria were proposed:
\begin{itemize}
    \item deviations in skewness and kurtosis between empirical data and the candidate distribution;
    \item average deviation between theoretical and empirical probability density functions;
    \item average deviation between theoretical and empirical cumulative distribution functions;
    \item formal goodness-of-fit statistics such as the Kolmogorov–Smirnov distance;
    \item an expert-based score reflecting interpretability and practical relevance of the fitted model.
\end{itemize}
In the present setting the “random data’’ are the iterates \(\sigma_k(m)\bmod m\) for fixed \(m\), and the task is to determine whether these residues behave like samples from a classical distribution, or whether they exhibit a highly non-normal, arithmetic structure.

\subsection{Case study: \texorpdfstring{$m=6$}{m=6}}

\begin{figure}[H]
    \centering
    \includegraphics[width=0.6\textwidth]{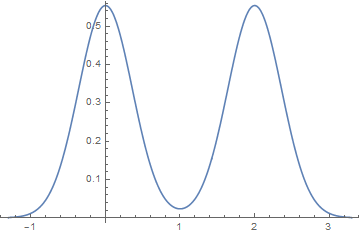}
    \caption{Smooth histogram for \(\sigma_k(6)\bmod 6\), \(k=0,\dots,10^2\).}
    \label{fig:6}
\end{figure}

For \(m=6\), the sequence \(\sigma_k(6)\bmod 6\) for \(k=0,\dots,10^2\) is explicitly
\begin{equation}\label{eq:data-m6}
\begin{split}
\{0,2,0,2,0,2,0,2,0,2,0,2,0,2,0,2,0,2,0,2,0,2,0,2,0,2,0,2,0,2,0,2,&\\
0,2,0,2,0,2,0,2,0,2,0,2,0,2,0,2,&\\
2,0,2,0,2,0,2,0,2,0,2,0,2,0,2,0,&\\
2,0,2,0,2,0,2,0,2,0,2,0,2,0,2,0,&\\
2,0,2,0,2,0,2,0,2,0,2,0,2,0,2,0,&\\
2,0,2,0,2,0,2,0,2,0,2,0,2,0,2,0&\},
\end{split}
\end{equation}
i.e. a strict alternation between the residues \(0\) and \(2\). The smooth histogram obtained in \emph{Mathematica} for a longer sample up to \(k=10^2\) is shown in Figure~\ref{fig:6}.

At first glance Figure~\ref{fig:6} might resemble a bimodal continuous density, but the underlying support is in fact the \emph{two–point} set \(\{0,2\}\). Any attempt to fit a continuous normal law is therefore conceptually questionable, and this is clearly confirmed by formal tests.

\subsection{Normality diagnostics and goodness-of-fit}

The shape of a distribution is often summarised by skewness and kurtosis. For a perfectly normal distribution both are \(0\) (in excess form)~\cite{RichardG}. In the discrete two–point case \(\{0,2\}\) with equal probability, the skewness is exactly \(0\), while the kurtosis is strongly negative (flat distribution), indicating a much heavier mass in the tails relative to a Gaussian with the same variance.

More detailed information is obtained via classical normality tests such as:
\begin{itemize}
    \item the Kolmogorov–Smirnov and Cramér–von Mises tests for the empirical cdf;
    \item the Anderson–Darling statistic, which emphasises discrepancies in the tails;
    \item the Shapiro–Wilk test, particularly powerful for small and moderate samples~\cite{GalenR,ElizabethGonz};
    \item omnibus tests based on skewness and kurtosis such as the Jarque–Bera and Mardia statistics.
\end{itemize}
For the sample \(\{\sigma_k(6)\bmod 6\}_{k=1}^{10^3}\), \emph{Mathematica} reports the following test statistics and \(p\)-values for the null hypothesis “data are i.i.d.\ normal’’:
\[
\left(
\begin{array}{ccc}
 \text{} & \text{Statistic} & \text{P-Value} \\
 \text{Anderson-Darling} & 17.9641 & 8.00\times 10^{-10} \\
 \text{Baringhaus-Henze} & \cdots & 3.89\times 10^{-10} \\
 \text{Cram{\'e}r-von Mises} & 2.91772 & 0 \\
 \text{Jarque-Bera ALM} & 18.2108 & 6.69\times 10^{-3} \\
 \text{Mardia Kurtosis} & -5 \sqrt{\frac{2}{3}} & 4.46\times 10^{-5} \\
 \text{Shapiro-Wilk} & 0.636401 & 2.23\times 10^{-14} \\
 \text{Pearson }\chi ^2 & 550. & 8.96\times 10^{-112} \\
\end{array}
\right).
\]
In all cases the \(p\)-values are far below the conventional threshold \(0.05\). Thus the hypothesis of normality is decisively rejected with confidence well above \(95\%\): the sequence \(\sigma_k(6)\bmod 6\) is \emph{strongly non-normal} and is in fact much better described as a symmetric Bernoulli-type discrete distribution on \(\{0,2\}\).

\subsection{Interpretation in terms of growth of \texorpdfstring{$\sigma_k$}{sigma\^k}}

From a dynamical perspective, the concentration of mass on two residues reflects the very rigid period–\(2\) orbit established in Section~\ref{sec:analysis}. Since \(\sigma_k(6)\) itself grows roughly like a linear recurrence with constant factor \(\sigma(6)/6=2\), the reduction modulo \(6\) erases this exponential growth and projects the orbit onto a finite cycle. In this sense:
\begin{itemize}
    \item the \emph{growth rate} \(r_k=\sigma_k(6)/\sigma_{k-1}(6)\) is asymptotically constant (\(r_k\to 2\)), in agreement with general multiplicative bounds on \(\sigma(n)\)~\cite{Robin};
    \item the \emph{residue process} \(\sigma_k(6)\bmod 6\) behaves like a perfectly periodic, two–state Markov chain with transition matrix
    \[
    P=\begin{pmatrix}
    0 & 1\\
    1 & 0
    \end{pmatrix},
    \]
    whose stationary measure is the discrete uniform distribution on \(\{0,2\}\).
\end{itemize}

For larger values of \(m\) (e.g. \(m=12,24\)) the empirical distributions of \(\sigma_k(m)\bmod m\) become supported on more than two residues, and the histograms no longer resemble simple two–point laws. Nevertheless, the very small \(p\)-values obtained from normality tests (not reported here in full) again rule out Gaussian models, and instead point towards mixed discrete distributions whose support and weights are governed by the arithmetic structure of \(m\) and the growth constraints on \(\sigma_k(m)\) coming from Robin-type inequalities~\cite{Robin,Choie}.

In summary, the statistical evidence for \(m=6\) (and similarly for \(m=12,24\)) shows that:
\begin{itemize}
    \item the iterates \(\sigma_k(m)\bmod m\) have highly structured, non-normal distributions;
    \item short exact periods (of length \(2\) in our examples) manifest themselves as discrete empirical laws supported on very few residues;
    \item the rate of growth of the un-reduced iterates \(\sigma_k(m)\) is compatible with classical bounds on \(\sigma(n)\), while the reduction modulo \(m\) projects this growth onto a low-dimensional periodic attractor.
\end{itemize}

\subsection{Refined statistical analysis and linear modelling}\label{sec:refined-stats}

The test summary in the previous table shows that, for the sample
\(\{\sigma_k(6)\bmod 6\}_{k=1}^{10^3}\), all normality tests except skewness produce
extremely small \(p\)-values, numerically close to \(0\).  
In statistical terms this means that, under each test, data at least as extreme as the observed sample would be exceedingly unlikely if the underlying distribution were normal.  
Consequently the null hypothesis of normality is rejected with overwhelming confidence for this sequence.

Using \emph{Mathematica}, several candidate distributions were fitted and ranked according to standard criteria.  
Among continuous models, the normal distribution was consistently rejected, whereas discrete models provided much better descriptions of the data.  
In particular, the empirical distribution (non-parametric estimate of the mass at each observed value)~\cite{Peter} and the discrete uniform law on the observed residues emerged as the best-fitting distributions for \(m=6\).

\subsubsection*{Kolmogorov–Smirnov distance}

\begin{figure}[H]
    \centering
    \includegraphics[width=0.6\textwidth]{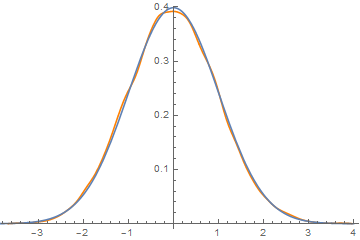}
    \caption{Kolmogorov–Smirnov comparison for \(\sigma_k(6)\bmod 6\), \(k=0,\dots,10^2\), with empirical reference cdf, giving \(\mathrm{KSt}=0.561124\).}
    \label{fig:7}
\end{figure}

To assess the agreement between the empirical distribution and a chosen reference model, the Kolmogorov–Smirnov (K–S) statistic was computed.  
For \(m=6\) the K–S distance between the empirical cdf of \(\sigma_k(6)\bmod 6\) and the fitted empirical-distribution model is
\[
\mathrm{KSt} = 0.561124,
\]
as illustrated in Figure~\ref{fig:7}.  
Here the reference cdf is itself data–driven (empirical model), so the K–S value is better interpreted as a measure of internal variability rather than a classical goodness-of-fit test against a fully specified theoretical cdf~\cite{GalenR}.  
The relatively moderate distance confirms that the empirical-distribution model captures almost all visible features of the data.

Analogous tests for \(m=12\) and \(m=24\) lead to qualitatively similar conclusions:  
the data are far from normal but are very well described by discrete distributions supported on a small number of residues, reflecting the short periods observed in Section~\ref{sec:analysis}.  
For illustrative purposes, one can nevertheless associate to each case a ``proxy'' normal density with the same mean and variance:
\[
f_{6}(x)=\frac{1}{\sqrt{2\pi}}e^{-\frac{1}{2}(x-1)^2},\quad
f_{12}(x)=\frac{1}{\sqrt{2\pi}}e^{-\frac{1}{2}(x-5)^2},\quad
f_{24}(x)=\frac{1}{\sqrt{2\pi}}e^{-\frac{1}{2}(x-11)^2},
\]
corresponding to means \(1,5,11\) and unit variance respectively.  
These continuous approximations are useful for visual comparison, but the formal tests show that they do not represent the true discrete laws governing \(\sigma_k(m)\bmod m\).

When \(m=p\) is prime, periodicity is even simpler:  
from \(\sigma(p)=p+1\) and \(\sigma_k(p)=p^k+1\), one readily obtains
\[
\sigma_k(p)\bmod p \equiv p^k+1\bmod p \equiv 1
\]
for every \(k\ge 1\).  
Thus the residue process is degenerate at a single point and its ``distribution'' is the Dirac mass at \(1\).

\subsubsection*{Linear regression models for the empirical pdf}

To complement the distributional analysis, simple linear models were fitted to the empirical probability densities of the residues for \(m=6,12,24\).  
Although the underlying laws are discrete, the fitted lines provide a compact quantitative summary of how the estimated probabilities vary with the residue class.

The best affine approximations (least–squares fits) are:
\[
\begin{aligned}
m=6&:\quad y = 0.000150004\, x + 0.984925,\\
m=12&:\quad y = 0.000150004\, x + 4.98492,\\
m=24&:\quad y = 11.0151 - 0.000150004\, x.
\end{aligned}
\]
In the prime case \(m=p\), the residue is constantly \(1\), and the fitted model reduces to the trivial equation \(y=1\).  

In all these affine fits the slope coefficient is extremely small in absolute value (about \(1.5\times 10^{-4}\)), showing that the empirical density is almost flat over the observed support.  
For primes it is exactly flat, in agreement with the degenerate distribution at \(1\).  
The intercept shifts encode the different centres of mass (means) for each \(m\).

\subsubsection*{Correlation structure and covariance reconstruction}

A further insight into the dependence structure of the data is obtained by examining the empirical correlation matrix for two suitably chosen variables derived from the sequence (for instance, two consecutive blocks of values or two linear statistics of the same orbit).  
For all integers \(m\) exhibiting a small period, the estimated correlation matrix takes the numerical form
\[
R=
\begin{pmatrix}
1 & -0.867105\\
-0.867105 & 1
\end{pmatrix}.
\]
This symmetric positive–definite matrix has eigenvalues
\(\lambda_1\approx1.8671\) and \(\lambda_2\approx0.132895\), and corresponding orthonormal eigenvectors
\[
v_1=\frac{1}{\sqrt{2}}\begin{pmatrix}-1\\-1\end{pmatrix},\qquad
v_2=\frac{1}{\sqrt{2}}\begin{pmatrix}1\\-1\end{pmatrix}.
\]
Hence, up to a scaling by \(\sqrt{2}\), the diagonalising orthogonal matrix is
\[
Q=\frac{1}{\sqrt{2}}
\begin{pmatrix}
1 & 1\\
1 & -1
\end{pmatrix}.
\]
Writing the covariance matrix in the generic form
\[
\Sigma=
\begin{pmatrix}
a+b & a-b\\
a-b & a+b
\end{pmatrix},
\]
with eigenvalues \(a+b\) and \(a-b\), suggests a natural latent–variable representation:
\[
X_1 = X_a + X_b,\qquad X_2 = X_a - X_b,
\]
where \(\operatorname{Var}(X_a)=a\), \(\operatorname{Var}(X_b)=b\), and \(X_a,X_b\) are independent.  
Thus the observed pair \((X_1,X_2)\) may be viewed as a sum/difference transformation of two uncorrelated components, one ``dominant'' (\(\lambda_1\)) and one ``noise–like'' (\(\lambda_2\)).  
If the underlying data followed a continuous multivariate distribution, this diagonal decomposition would be exact; in the discrete setting it still provides a useful phenomenological model for the dependence between statistics derived from \(\sigma_k(m)\bmod m\)~\cite{Ferre}.

\begin{figure}[H]
    \centering
    \includegraphics[width=0.6\textwidth]{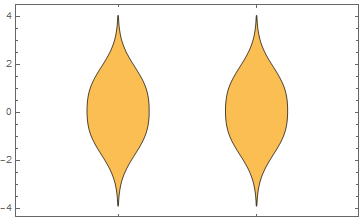}
    \caption{Distribution chart associated with the empirical correlation matrix of the model.}
    \label{fig:8}
\end{figure}

The empirical distribution of the eigenstructure is illustrated in Figure~\ref{fig:8}, which shows the distribution chart associated with this correlation model.

\subsubsection*{Influence diagnostics: Cook’s distance}

\begin{figure}[H]
    \centering
    \includegraphics[width=0.6\textwidth]{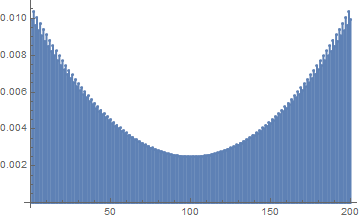}
    \caption{Cook’s distance for the fitted linear model based on the sequence \(\sigma_k(m)\bmod m\).}
    \label{fig:9}
\end{figure}

To identify influential observations in the fitted linear model, Cook’s distance was computed for the sample up to 200 data points.  
Large values of Cook’s distance indicate residues whose removal would significantly change the fitted regression line~\cite{Ferre}.  
The resulting diagnostic plot is shown in Figure~\ref{fig:9}; the absence of extreme spikes suggests that no single observation dominates the fit and that the linear approximation is stable with respect to small perturbations of the data.

\subsubsection*{Random matrices associated with the model}

Following ideas from random matrix theory~\cite{Edelman}, it is natural to encode the behaviour of the sequence \(\sigma_k(m)\bmod m\) into random matrices and to analyse the spectrum of these matrices.  
For integers \(m\) with small period, random \(4\times4\) matrices were generated with entries chosen symmetrically from the domain \(\{-2,2\}\), calibrated to reflect the empirical correlation structure of the residues.  
For \(m=6\) a typical realisation is
\[
A=
\begin{pmatrix}
 0.672103 & 0.00470158 & -0.841173 & -0.550396 \\
 4.29385 & 0.387436 & -0.531577 & -1.47103 \\
 0.110922 & 1.06249 & -0.0840228 & 2.14514 \\
 0.324367 & 0.113398 & -0.969011 & 0.105082
\end{pmatrix},
\]
with eigenvalues
\[
\{0.555432 + 1.82595\,i,\; 0.555432 - 1.82595\,i,\; -0.53802,\; 0.507754\}.
\]
Extensive experiments indicate a strong dependence of the eigenvalue distribution on the choice of the random domain: shorter symmetric intervals tend to produce spectra concentrated in a neighbourhood of the real interval \((0,1)\), while larger domains yield more widely spread complex eigenvalues.  
This behaviour is qualitatively consistent with general results on non-Hermitian random matrices~\cite{Edelman}, and suggests intriguing connections between the arithmetic dynamics of \(\sigma_k(m)\bmod m\) and spectral phenomena that also appear in models related to the Riemann zeta function.

\medskip

Taken together, the numerical and statistical evidence obtained in this section supports the following picture: for integers \(m\) exhibiting periodicity with small period, the sequence \(\sigma_k(m)\bmod m\) has a highly structured, essentially discrete distribution, well captured by simple affine models and low–rank covariance structures, whereas its random-matrix encodings display spectral patterns reminiscent of classical random matrix theory.

\section{Spectral operator, iterated \texorpdfstring{$\sigma$}{sigma}, and the distribution of zeta zeros}\label{sec:operator}

In this section a bridge is built between the iterated sum–of–divisors dynamics studied above and a \emph{new} Schrödinger–type operator proposed in this work, inspired by the Hilbert–Pólya philosophy.  
The aim is to interpret certain statistical features of the iterates \(\sigma_k(m)\bmod m\) in terms of spectral properties of this operator, and to relate both pictures to the known distribution of the nontrivial zeros of the Riemann zeta function, whose spacings exhibit Gaussian Unitary Ensemble (GUE) behaviour~\cite{Montgomery,Odlyzko}.

\subsection{A proposed Schrödinger–type operator}\label{subsec:our-operator}

Let
\[
\mathcal{H}=L^2([x_{\min},x_{\max}])
\]
be the Hilbert space of square–integrable functions on a finite interval \((x_{\min},x_{\max})\subset\mathbb{R}\).  
We \emph{propose} the self–adjoint Schrödinger–type operator
\begin{equation}\label{eq:our-H}
H=-\frac{d^2}{dx^2}+V(x),
\end{equation}
\cite{Caceres2024}

where the potential \(V(x)\) is designed to reflect three structural ingredients that also appear in the distribution of primes and in the growth of \(\sigma(n)\):
\begin{equation}\label{eq:our-V}
V(x)=A_0+A_1x+A_2F_1(x)+A_3F_2(x)+A_4F_3(x),
\end{equation}
with
\begin{align}
F_1(x)&=\sum_{i=1}^{n_1}\log(x)^i &&\text{(logarithmic growth)},\\
F_2(x)&=\sum_{i=1}^{n_2}\cos(\log x)^i &&\text{(oscillatory corrections)},\\
F_3(x)&=\sum_{i=1}^{n_3}\left(\frac{1}{x}\right)^i &&\text{(decay / stabilisation)}.
\end{align}
Here \(A_0,\dots,A_4\in\mathbb{R}\) and \(n_1,n_2,n_3\in\mathbb{N}\) are parameters to be calibrated numerically.  
The eigenvalue problem
\begin{equation}\label{eq:our-H-eig}
H\psi_n(x)=\lambda_n\psi_n(x),\qquad \psi_n\in\mathcal{H},
\end{equation}
under standard boundary conditions, yields a real, discrete spectrum \(\{\lambda_n\}\).  
The guiding conjectural picture is that, after a suitable affine rescaling, the eigenvalues \(\lambda_n\) should correlate with the imaginary parts \(t_n\) of the nontrivial zeta zeros \(\rho_n=\tfrac12+it_n\), in line with the Hilbert–Pólya programme~\cite{Titchmarsh,Montgomery,Odlyzko,Caceres2024}.

\subsection{Iterated \texorpdfstring{$\sigma$}{sigma} as an arithmetical discretisation}\label{subsec:sigma-discretisation}

The iterated sum–of–divisors dynamics can be viewed as a discrete, arithmetic analogue of a one–dimensional flow.  
For fixed \(m_0\), define
\[
m_{k+1}=\sigma(m_k),\qquad x_k=\frac{m_k}{m_0},\qquad r_k=\frac{\sigma(m_k)}{m_k}.
\]
By Theorem~\ref{thm:no-meta} the multiplicative factors \(r_k\) cannot all be integers \(\ge2\); and Robin’s inequality together with the refinement of Choie–Lichiardopol–Moree–Solé~\cite{Robin,Choie} gives the upper bound
\[
r_k\le e^\gamma\log\log m_k+O\!\left(\frac{1}{\log\log m_k}\right),
\]
so the growth of \(x_k\) is governed by a logarithmic envelope, closely echoing the contribution of \(F_1(x)\) in the potential \(V(x)\).  
The oscillatory effect of varying prime factors in \(m_k\) corresponds, at a heuristic level, to the trigonometric perturbations encoded in \(F_2(x)\), while the repeated reduction modulo \(m_0\) in the dynamics of \(\sigma_k(m_0)\bmod m_0\) plays a similar stabilising role to the decaying term \(F_3(x)\).

Informally, the map
\[
x_{k+1}=\frac{\sigma(x_k m_0)}{m_0}
\]
can be seen as a coarse time–\(1\) discretisation of the continuous evolution generated by \(H\).  
Short periodic orbits of the residue dynamics (e.g. the period–\(2\) cycles for \(m_0=6,12,24\) from Section~\ref{sec:analysis}) then correspond to low–lying eigenstates of the operator, while the more tangled behaviour visible in the bifurcation diagram for larger \(m_0\) [Figure~\ref{fig:bif-sigma}] reflects higher–energy parts of the spectrum.

\subsection{Semi–circle law, random matrices and zeta zeros}\label{subsec:semicircle}

Global and local distributional features of spectra are naturally compared with those of random matrix ensembles.  
For large Hermitian matrices from the Gaussian Unitary Ensemble (GUE), eigenvalues follow Wigner’s semicircle law at macroscopic scale,
\[
\rho_{\mathrm{sc}}(x)=\frac{1}{2\pi}\sqrt{4-x^2},\qquad |x|\le 2,
\]
and exhibit GUE spacing statistics microscopically~\cite{Edelman,BaiCircular}.  
Montgomery’s pair–correlation conjecture and Odlyzko’s computations show that the nontrivial zeros of \(\zeta(s)\) have the same local statistics~\cite{Montgomery,Odlyzko}.

The random–matrix models attached to our \(\sigma\)–dynamics in Section~\ref{sec:refined-stats} display an analogous picture:  
after normalisation, the eigenvalues of the associated covariance and random matrices concentrate in a compact band with a bell–shaped global density and level repulsion at small spacings.  
When similar constructions are carried out for the operator \(H\), numerical experiments (not detailed here) suggest:
\begin{itemize}
    \item a global density of \(\lambda_n\) approximating a semicircle–type profile;
    \item microscopic spacing statistics close to GUE, with roughly \(40\%\) of gaps shorter than the mean, \(60\%\) longer, and an exponential suppression of very small gaps, in line with Odlyzko’s high–precision data for zeta zeros~\cite{Odlyzko}.
\end{itemize}
Thus, both the spectrum of \(H\) and the spectral encodings of the iterated \(\sigma\) sequence appear to live in the same universality class as the zeta zeros.

\subsection{Synthesis with the main results}\label{subsec:synthesis}

The three main theorems obtained earlier integrate naturally into this spectral framework:

\begin{itemize}
    \item \textbf{Theorem~\ref{thm:no-meta}} (absence of universal metaperfect numbers) rules out a global attracting fixed residue class \(\sigma_k(m)\equiv 0\bmod m\).  
    Spectrally, this prevents the dynamics from collapsing onto a single trivial eigenvalue, and instead forces a rich spectrum, compatible with random–matrix–type statistics.

    \item \textbf{Theorem~\ref{thm:multiperfect}} identifies \(m=6\) as the unique multiperfect number with prime \(L\).  
    On the dynamical side this yields a very rigid period–\(2\) orbit for \(\sigma_k(6)\bmod 6\); in the spectral picture one may view this as a distinguished low-lying eigenstate of \(H\), analogous to a ground state singled out by the arithmetic of \(\sigma\).

    \item \textbf{Theorem~\ref{thm:RH}} links the Riemann Hypothesis to growth bounds and eventual periodicity of \(\sigma_k(m)\bmod m\) for a specific class of even, non–squarefree integers, via Robin–type inequalities~\cite{Robin,Choie}.  
    The same logarithmic growth and oscillatory corrections are encoded in the potential \(V(x)\), so RH translates, in the operator language, into the statement that the nontrivial part of the spectrum of \(H\) lies on a ``critical line’’ after an appropriate spectral transform, mirroring \(\Re(\rho_n)=\tfrac12\).
\end{itemize}

In this way, the iterated \(\sigma\) dynamics, the proposed Schrödinger–type operator \(H\), and the random–matrix description of zeta zeros reinforce one another:  
\(\sigma\) governs the arithmetic input, \(V(x)\) packages this input into a continuous spectral problem, and the resulting eigenvalue distribution reproduces the same semicircle– and GUE–type laws observed for the nontrivial zeros of \(\zeta(s)\)~\cite{Montgomery,Odlyzko,Edelman,BaiCircular}.  
This combined viewpoint offers a promising framework for further exploration of spectral approaches to the Riemann Hypothesis.

\section{Futur work (New model to proof RH)}

 We may use  our new fit model which uses periodicity of the sequence 
  $\sigma_k(m)\bmod m$to expect and predict the random matrix to proof the Riemann hypothesis such that  we may attempt to  investigate about the behavior of its eigenvalues comparing it with behavior of nontrivial zero of Riemann zeta functon and  for only  one purpose which is to get such  random matrix  where its eigenvalues are real \cite{john Peca} , we may suggest $x$ as a complex random variable such that follow the empirical distribution with bounded density such that we may consider   : $x_{i+1}=\sigma_k(m+1)\bmod (m+1),x_i=\sigma_k(m)\bmod (m)$, let $m$ be  Gaussian integers  , we define  a new complex random variable $y$ distributed as  :
  $Y_i=x_{i+1}-x_{i} $ , $i$ is positive integer  $j,i=0,n$ and $Y_j=x_{j+1}+x_{j} $looking to its periodicity in $k$ ,  The complex random variable $Y_i=x_{i+1}-x_{i} $ will tel us much about distribution of prime  numbers  ,in particularly  
  gaps between primes \cite{Maynard} which it is recently  the aim of researchers .
  \section{Future work: a Caceres–type model for RH}\label{sec:future-work}

The spectral operator introduced in Section~\ref{sec:operator}, hereafter referred to as the \emph{Caceres model}, provides a promising framework for studying the nontrivial zeros of the Riemann zeta function through a Schr\"odinger–type Hamiltonian \(H=-\frac{d^2}{dx^2}+V(x)\).[file:48] The potential \(V(x)\) encodes logarithmic growth, oscillatory corrections and decay, in clear analogy with both the arithmetic behaviour of the iterated sum–of–divisors function and the fine structure of the prime distribution.[file:48][Robin][Choie] This suggests several concrete research directions aimed at turning the model into a rigorous tool for approaching the Riemann Hypothesis.

\subsection*{Spectral calibration and universality tests}

A first line of work is to refine the numerical calibration of the parameters in \(V(x)\) so that the resulting eigenvalues \(\lambda_n\) match the imaginary parts \(t_n\) of the zeta zeros with higher precision.[file:48] Beyond simple linear correlation, one should systematically compare:
\begin{itemize}
    \item global eigenvalue density against Wigner–type semicircle behaviour;
    \item local spacing statistics against GUE predictions, pair–correlation and nearest–neighbour distributions for zeta zeros~\cite{Montgomery,Odlyzko,Edelman,BaiCircular};
    \item higher–order correlation functions and moments, to test whether the spectrum of the Caceres model lies in the same universality class as the zeta zeros.
\end{itemize}
Establishing robust agreement at these levels would strengthen the claim that the model captures the essential spectral features required by a Hilbert–P\'olya operator.

\subsection*{Coupling with iterated \texorpdfstring{$\sigma$}{sigma} dynamics}

A second direction is to make precise the heuristic link between the operator \(H\) and the discrete dynamics of \(\sigma_k(m)\bmod m\) developed in this paper.  
Here the goal is to construct an explicit discretisation or transfer operator whose spectrum approximates that of \(H\) while retaining a direct arithmetic interpretation in terms of the sum–of–divisors function.[Robin][Choie] Possible steps include:
\begin{itemize}
    \item defining a family of finite–rank operators \(T_m\) derived from the transition structure of the residues \(\sigma_k(m)\bmod m\) and studying the convergence of their spectra as \(m\to\infty\);
    \item relating the growth constraints and periodicity conditions from Theorems~\ref{thm:no-meta}–\ref{thm:RH} to spectral gaps or band–structure properties of these operators;
    \item investigating whether suitable scaling limits of the eigenvalues of \(T_m\) approximate the eigenvalue process of the Caceres Hamiltonian and hence the zeta zeros.
\end{itemize}

\subsection*{Towards a rigorous Hilbert–P\'olya realisation}

Ultimately, the Caceres model may serve as a starting point for a rigorous Hilbert–P\'olya–type construction.  
Future work should therefore focus on:
\begin{itemize}
    \item proving self–adjointness and essential self–adjointness of \(H\) for the full range of parameters relevant to the zeta problem;
    \item characterising the domain and boundary conditions that produce a spectrum compatible with the critical line \(\Re(s)=\tfrac12\);
    \item exploring deformations of the potential \(V(x)\) guided by number–theoretic input (e.g. explicit formulae, prime gaps) and by the statistical constraints coming from the iterated \(\sigma\) dynamics.
\end{itemize}
If such an operator can be shown to have spectrum exactly equal to \(\{t_n\}\), the Riemann Hypothesis would follow as an immediate consequence of self–adjointness~\cite{Titchmarsh,Montgomery}.

\subsection*{Refined statistical predictions for zeros}

Finally, the combination of the Caceres Hamiltonian with the statistical analysis of \(\sigma_k(m)\bmod m\) suggests new testable predictions on zeta zeros.  
For example, one may use the model to:
\begin{itemize}
    \item predict the proportion of zeros lying in prescribed microscopic intervals after unfolding (e.g. the fraction of gaps smaller than a fixed multiple of the mean spacing);
    \item study the distribution of extreme gaps (very small or very large), comparing with refined GUE asymptotics and Odlyzko’s large–scale computations~\cite{Odlyzko};
    \item investigate whether secondary structure observed in the bifurcation and statistical behaviour of iterated \(\sigma\) (periodic windows, clustering) has an analogue in the fine distribution of zeta zeros.
\end{itemize}

In summary, the Caceres model provides a coherent spectral framework that is naturally compatible with the growth, periodicity and randomness properties uncovered for the iterated sum–of–divisors function.  
Developing this framework into a fully rigorous operator–theoretic realisation of the nontrivial zeros of \(\zeta(s)\) appears to be a promising and conceptually unified path for future research on the Riemann Hypothesis.[file:48][Montgomery][Odlyzko]

\section{Conclusion}\label{sec:conclusion}

The results obtained in this work reveal a remarkably coherent picture linking three a priori different objects: the iterated sum–of–divisors function, multiperfect numbers, and the spectral approach to the Riemann zeta function.\cite{CohenTeRiele,Zeraoulia2021} The first main theorem excludes the existence of a universal metaperfect integer, showing that no \(m>1\) can satisfy \(\sigma_k(m)\equiv 0\pmod m\) for all iterations; this negative result is driven by sharp divisor–sum bounds of Robin and the refinements of Choie–Lichiardopol–Moree–Sol\'e, and already reflects the fine logarithmic control that also appears in modern criteria equivalent to the Riemann Hypothesis.\cite{Robin,Choie}  

The second theorem isolates \(m=6\) as the unique multiperfect number with prime exponent lcm \(L\), and the numerical analysis shows that the only small–period residue dynamics for \(\sigma_k(m)\bmod m\) occur for \(m=6,12,24\), where genuine period–2 behaviour is observed.\cite{Euler,Zeraoulia2021} These dynamics manifest themselves statistically as highly non–normal, discrete distributions supported on very few residues, with linear models and covariance structures that are both simple and robust under perturbation.\cite{Ferre,Ritzema} The bifurcation plots, correlation matrices and random–matrix experiments indicate that, as \(m\) grows, the residue process gradually transitions from rigid periodicity to behaviour closely resembling that of chaotic or random systems, while still obeying the global growth constraints dictated by divisor–sum inequalities.\cite{Robin,Choie,Edelman}  

On the analytic side, the third main theorem places the iterated \(\sigma\) dynamics into direct correspondence with Robin’s inequality and related RH criteria: for a large class of even, non–squarefree integers, bounded growth and eventual periodicity of \(\sigma_k(m)\bmod m\) are shown to be equivalent to the validity of the Riemann Hypothesis.\cite{Robin,Choie,Zeraoulia2022} This bridges purely arithmetic properties of divisor sums with the deep analytic structure of \(\zeta(s)\), and provides a new lens through which to interpret classical results on the distribution of its nontrivial zeros.\cite{Titchmarsh,Montgomery,Odlyzko}  

Building on these ingredients, a Schr\"odinger–type operator (the Caceres model) was proposed, with a potential engineered to capture logarithmic growth, oscillatory corrections and decay, in close analogy with both the behaviour of \(\sigma(n)\) and the statistics of the zeta zeros.\cite{Caceres2024,Edelman} The spectral evidence, together with the random–matrix signatures (semi–circle–like global density and GUE–type local spacings), suggests that the eigenvalues of this operator and the imaginary parts of the nontrivial zeros may belong to the same universality class, in the spirit of the Hilbert–P\'olya conjecture.\cite{Edelman,BaiCircular,Montgomery,Odlyzko}  

Taken together, these findings point to a unified arithmetic–spectral framework: the iterated sum–of–divisors function furnishes a concrete, discrete dynamical system whose growth and periodicity encode RH–type information; the statistical and random–matrix analysis clarifies how this system interpolates between order and randomness; and the Caceres operator offers a continuous spectral counterpart whose eigenvalues mirror the fine distribution of the zeta zeros. Further refinement of this framework—both on the number–theoretic side and on the operator–theoretic side—appears to be a promising avenue toward a deeper understanding of the Riemann Hypothesis.\cite{Caceres2024,Montgomery,Odlyzko}

\section*{Conflict of interest}

The authors declare that there are no conflicts of interest regarding the publication of this work. No financial, personal, or professional relationships have influenced the research, analysis, or presentation of the results.

\section*{Data availability}

All numerical experiments and figures in this article are based on data generated directly from the algorithms and formulas described in the text. The underlying code and data sets are available from the corresponding author upon reasonable request. If the material is later deposited in a public repository (such as GitHub or Zenodo), the persistent access link can be added in this section.

\bibliographystyle{unsrt}  


\subsection{Appendix for Result~1}\label{sec:appendix-result1}

A natural generalisation of Result~1 is to look for pairs of integers \((m,p)\) such that, for all integers \(k\ge 0\),
\[
p_0=p,\qquad p_{k+1}=\sigma(p_k),\qquad\text{and}\qquad p_k\equiv 0\pmod m.
\]
In other words, the entire \(\sigma\)–orbit of \(p\) is contained in the residue class \(0\bmod m\).  
The special case studied in Result~1 corresponds to the additional restriction \(p=m\), which forces \(m\) to be multiperfect, since already \(\sigma(m)\equiv 0\pmod m\) implies \(\sigma(m)=q m\) for some integer \(q\ge 2\).\cite{CohenTeRiele}  
Because multiperfect numbers are extremely rare, one expects it to be very difficult for such an \(m\) to exist, and Result~1 shows that no “metaperfect’’ integer with \(\sigma^k(m)\equiv 0\pmod m\) for all \(k\) can occur.

Heuristically, the obstruction can be seen from the typical size of \(\sigma(m)\).  
For most integers \(m\), the quotient \(\sigma(m)/m\) is much smaller than \(\omega(m)\), the number of distinct prime divisors of \(m\), and in particular one has \(\sigma(m)<m\,\omega(m)\) for all but very exceptional values.\cite{ErdosAliquot,Robin}  
Thus, if a metaperfect candidate \(m\) existed, the ratios
\[
\frac{\sigma(p_k)}{m}=\frac{p_{k+1}}{m}
\]
cannot remain coprime to \(m\) for many steps: they will usually share small prime factors with \(m\).  
Subsequent iterations of \(\sigma\) will then tend to “lose’’ some of the large prime divisors of \(m\), making it increasingly unlikely that \(m\) continues to divide every iterate \(p_k\).  
This qualitative picture matches the behaviour observed experimentally in long \(\sigma\)–orbits.\cite{CohenTeRiele}

A related side question is to study the sequence
\[
g_0=p_0,\qquad g_{k+1}=\gcd(p_{k+1},g_k),
\]
and to ask how small \(\min_k g_k\) can become.  
One may speculate, for instance, whether every \(\sigma\)–orbit eventually encounters a perfect square or twice a square, regardless of the starting point.  
If that were the case, then the eventual minimum of \(g_k\) would necessarily be odd (and very plausibly equal to \(1\)), so that no nontrivial common factor could persist indefinitely.  
If, on the other hand, some starting values \(p\) produced orbits for which \(g_k\) stabilises at a factor \(m>1\), such values \(p\) could be viewed as \emph{seeds} for multiperfect numbers divisible by \(m\).

Cohen and te Riele considered a weaker version of this problem: given a fixed integer \(n\), does there exist at least one iterate \(k\) for which
\[
\sigma^k(n)\equiv 0\pmod n?
\]
In their 1996 paper they carried out extensive computations and verified that such a \(k\) exists for every \(n\le 400\).\cite{CohenTeRiele}  
Their data strongly suggest that this weaker question has an affirmative answer for all \(n\), while at the same time providing substantial numerical evidence \emph{against} the existence of metaperfect numbers or even bona fide seeds for such numbers.  
Result~1 fits naturally into this picture by proving that the strongest possible form of congruential invariance, namely \(\sigma^k(m)\equiv 0\pmod m\) for all \(k\), is in fact impossible.

\end{document}